\documentclass[11pt,a4paper]{article}

\usepackage[top=2.4cm,bottom=2.4cm,left=2.6cm,right=2.6cm]{geometry}
\usepackage{amsmath,amssymb,amsthm}
\usepackage{mathtools}
\usepackage{authblk}
\usepackage[giveninits=true,maxbibnames=10,url=false,isbn=false]{biblatex}
\usepackage[hidelinks]{hyperref}

\newcommand{\dd}{\mathrm{d}}
\newcommand{\Var}{\operatorname{Var}}

\newcommand{\Corr}{\operatorname{Corr}}
\newcommand{\Ai}{\operatorname{Ai}}
\newcommand{\Bi}{\operatorname{Bi}}
\newcommand{\E}{\mathbb{E}}

\DeclareMathOperator*{\argmin}{\arg\,min}

\title{Forecasting tipping times from data: extrapolation of critical transitions and uncertainty quantification}
\author[1]{Davide Papapicco}
\author[1]{Lauren Smith}
\author[1]{Graham Donovan}
\affil[1]{Department of Mathematics, University of Auckland, New Zealand}
\date{}

\begin{document}

\maketitle

\begin{abstract}\noindent
        The long--term behaviour of complex systems in ecology, biology and climate science can sometimes experience nonlinear excursions away from a regime of stable evolution around equilibria.
        These events, often called critical transitions or tipping phenomena, are usually associated with the loss of resilience of a given system to recover from external perturbations.
        Due to the potentially catastrophic repercussions that critical transitions can bring upon the sustainability of complex systems, a great deal of effort has recently being devoted in the characterization and early detection of tipping events from timeseries.
        In this paper we derive a fully data--driven pipeline for the accurate forecast of tipping times in complex systems.
        We combine a range of existing techniques in nonautonomous dynamics and statistical estimation for stochastic differential equations in a novel extrapolation scheme while providing accuracy and uncertainty quantification.
        Our aim is to provide an out--of--the--box methodology for prediction of critical transitions as well as the generation and collection of data to scientists and modelers.
        This pipeline relies on a series of assumptions, which we discuss throughout this paper, that provide an optimal observation protocol for the forecast of tipping times.
\end{abstract}

\section{Introduction}\label{sec:intro}
Tipping phenomena (or critical transitions) are abrupt changes in the state of complex systems that correspond to disruptive and potentially permanent alterations in their long--term behaviour.
They have been studied across a broad and diverse collection of natural and human critical phenomena including ecosystems degradation \cite{OSullivan2023,Kumar2026}, neurobiological activity spikes \cite{Boerner2024}, onset of asthma \cite{Venegas2005,Donovan2015}, weakening of ocean currents \cite{Boers2021,Dijkstra2026}, hysteresis in laser systems \cite{Hohl1995}, generation of instabilities in magnetically confined plasma \cite{Cowley2015} and thermoacustic systems \cite{Pavithran2021}, marine vessels capsizing \cite{Shao2026}, formation of traffic congenstion \cite{Chattopadhyay2026} and civilization collapse \cite{Scheffer2021,Doak2026}.

Catastrophic events are usually characterized by the loss of stability of steady--states due to nonlinear and noisy responses to slowly--varying external forcing.
Their canonical classification, proposed in \cite{Ashwin2012}, distinguishes between bifurcation--induced (B-), rate--induced (R-) and noise--induced (N-) tipping points.
Furthermore, slow--fast dynamics and bifurcation theory provide a modeling framework for the forecast of critical transitions \cite{Kuehn2011}.
The prediction of approaching catastrophes, known as \emph{early-warning signals}, usually relies on the detection of characteristic precursors from non--stationary timeseries.
The most common of these signatures is known as \emph{critical slowing down} \cite{Tredicce2004} and it is often associated with an increase in the variance and autocorrelation of a timeseries sample.

Many early-warning signals have been proposed in the last two decades; we refer interested readers to \cite{Boettiger2013}, and more recently \cite{Dakos2024}, for a concise overview on the main research on the topic.
Debates regarding the legitimacy of theoretical early-warnings applied to real--world data have since followed \cite{Ditlevsen2010,Bos2014,Dakos2015,Wilkat2019}, suggesting that caution should be exercised when relying on these indicators to predict critical events \cite{Morr2024a,Cui2025,Parry2026,Ditlevsen2026}.
In particular, pertinent concerns to applications are:
\begin{itemize}
  \item recording optimal observables from high--dimensional (spatio--temporal or networked) systems is necessary in order to detect critical slowing down \cite{Lohmann2025b};
  \item reliable and skillfull early-warnings require sufficient timescale separation between the fast (state) and the slow (forcing) variables \cite{Kuehn2011};
  \item predictions based on ensemble averages without uncertainty quantification of the single realizations can lead to missed warnings (i.e. false negatives) and/or spurious signals (i.e. false positives) \cite{Ashwin2025};
  \item lack of knowledge of baseline values of classical early-warning signals result in the lack of interpretability of statistical quanties without context and historical trends \cite{Papapicco2026}.
\end{itemize}

In this manuscript we are concerned with the problem of forecasting the time at which complex systems undergo critical transitions given past and present data.
The forecast is a (linear) extrapolation of the estimated return rate of nonlinear systems that are slowly ramped towards a saddle-node bifurcation.
Such saddle-node B-tipping events are found in many of the applications mentioned above in which the collected timeseries samples provide a discrete time approximation of a slow passage through the bifurcation \cite{Tsai2025}.
We provide a data--driven pipeline and an observation protocol with which scientists and modelers can collect data to predict such tipping events, far enough in the future and with sufficient accuracy, to enable policies that stop or reverse the catastrophe.

To achive this we leverage well--known results in classical and nonautonomous dynamical systems theory and stochastic dynamics as well as canonical statistical estimation techniques applied to timeseries data.
The pipeline and the protocol rely on a set of assumptions that the data must satisfy in order for the forecasting technique to yield consistent and reliable results.
Specifically, we discuss how the level of noise, rate of change of the forcing parameter, distance from the bifurcation and sample length can be optimally selected to robustly forecast tipping times.

The paper is structured as follows.
In \S~\ref{sec:tipping} we introduce the mathematical formalism of tipping times in slow--fast stochastic processes.
In \S~\ref{sec:formulation} we formulate the extrapolation problem at the core of the data--driven pipeline by presenting three techniques to estimate the return rate of a given sample, two of which have already been proposed in the literature and a new one which we introduce.
In \S~\ref{sec:observation} we derive the observation protocol, that is a series of conditions that allow for the accurate forecast of tipping times from collected measurements or experimental data.
In \S~\ref{sec:application} we showcase the performance of the pipeline on both low-- and high--dimensional timeseries.
In \S~\ref{sec:conclusion} we discuss our findings and suggest future directions for this line of research.

\section{Tipping times of slowly ramped systems}\label{sec:tipping}
The prototypical case of critical transitions in complex system is that of a dynamic saddle-node normal form, with slow (linear) passage through the bifurcation and (additive) noisy fluctuations modelled as an It\^o stochastic differential equation of the form
\begin{equation}\label{eq:model}
\dd X_t=f\bigl(X_t,\mu(t)\bigr)\,\dd t+\sigma\,\dd W_t,
\qquad
f(x,\mu):=-\mu-x^2,
\qquad
\mu(t):=\mu_0+\varepsilon t ,
\end{equation}
where $X_t\in\mathbb R$ is the state of the system, perturbed by a Wiener process $W_t$ with noise amplitude $\sigma>0$, and $\mu(t)$ is the bifurcation parameter, ramped with constant speed $\varepsilon>0$ from an initial value $\mu_0<0$.
The drift term $f(x,\mu)$ in \eqref{eq:model} has two equilibria $x_{\pm}(\mu)=\pm\sqrt{-\mu}$ for any $\mu<0$; at $\mu=0$ the two branches collide $x_+(0)=x_-(0)=0$ in a saddle-node bifurcation and for $\mu>0$ they annihilate each other.
Small perturbations around the stable equilibrium $x_+(\mu)$ decay exponentially fast at a rate dictated by the leading eigenvalue of the linearized vector field
\begin{equation}\label{eq:eigenvalue}
        \lambda(\mu)\coloneqq\partial_xf(x_+(\mu),\mu)=-2\sqrt{-\mu}<0 , \quad \mu<0 .
\end{equation}
Specifically, the quantity $1/|\lambda|$ identifies the \emph{correlation} or \emph{relaxation time} of noisy fluctuations towards $x_+(\mu)$.
Evidently, $|\lambda|\to0$ as $\mu\to0$, which in turn gives asymptotically longer relaxation time and eventually the blow-up $1/|\lambda|\to\infty$. 
This critical slowing down is a generic behaviour of any complex system driven close enough to a saddle-node bifurcation as the nonlinear dynamics is reduced to the normal form \eqref{eq:model} by a smooth change in coordinates.

For a given ramp speed $\varepsilon>0$, the time at which the leading eigenvalue \eqref{eq:eigenvalue} vanishes corresponds to the time at which the forcing parameter $\mu(t)$ reaches the saddle-node bifurcation $\mu(t)=0$, which is $\tau\coloneqq-\mu_0/\varepsilon>0$.
One can thus rewrite the linear ramp of the forcing parameter to be independent of the initial condition $\mu_0<0$ and parametrized by the time horizon $\tau$ at which the saddle-node bifurcation is reached, i.e.
\begin{equation}\label{eq:breaking_time}
        \mu(t) = \mu_0 + \varepsilon t = \varepsilon(t-\tau)<0\quad\forall t<\tau .
\end{equation}
Substituting \eqref{eq:breaking_time} in \eqref{eq:eigenvalue} gives an explicit time dependence of the approach to the bifurcation
\begin{equation}\label{eq:eigenvalue_time}
        \lambda(t) = -2\sqrt{\varepsilon(\tau-t)}\to0,\quad t\to\tau .
\end{equation}

\subsection{Deterministic tipping time}\label{subsec:delay}
Let us first consider the deterministic limit $\sigma\to0$ for which the nonautonomous system \eqref{eq:model} is known as the \emph{dynamic} saddle-node normal form \cite{Berglund2006}.
The time horizon $\tau$ (henceforth \emph{breaking time}) identifies the loss of hyperbolicity of the stable equilibrium branch $x_+(t)\coloneqq x_+(\mu(t))$ as per \eqref{eq:eigenvalue_time}.
In the nonautonomous regime however, the classic notion of $x_+(t)$ is replaced by the pullback attractor $x_{\mathrm{pb}}(t)$ which is the unique solution of \eqref{eq:model} (in the absence of noise $\sigma=0$) satisfying 
\begin{equation}\label{eq:pb_attractor}
  \lim_{t\to-\infty}|x_{\mathrm{pb}}(t)-x_+(t)|=0 .
\end{equation}
For a slow passage through the saddle-node bifurcation, i.e. $0<\varepsilon\ll1$, provided that $\mu_0=\lim_{t\to-\infty}\mu(t)\ll0$, the pullback attractor converges to a slow manifold \cite[Thm.~$2.1.7$]{Berglund2006}, which tracks the stable equilibrium branch $x_+(t)$ and persists past the breaking time $\tau$ for some time before diverging to $-\infty$.
The trajectory of a slow passage through a saddle-node bifurcation will thus undergo B-tipping at a delayed time $\tau_\mathrm{pb}>\tau$ (henceforth \emph{tipping time}), creating a window of opportunity $\Delta\tau\coloneqq\tau_{\mathrm{pb}}-\tau>0$ to reverse the forcing $\varepsilon$ and avoid the critical transition.

To quantify the tipping time, and thus identify the window of opportunity, we follow the procedure outlined in \cite{Li2019}.
We start by writing the dynamic saddle-node normal form as a Riccati equation
\begin{equation}\label{eq:riccati}
        \frac{\dd x}{\dd t} = f(x,\mu(t)) = \sum_{k=0}^{2}c_k(t)x^k,
\end{equation}
with $c_0(t)=\varepsilon(\tau - t)$, $c_1(t)=0$ and $c_2(t)=-1$, and convert it into a second-order, linear system via change of variables $x\mapsto -(c_2(t)y)^{-1}(\dd y/\dd t)$
\begin{equation}\label{eq:converted}
        \frac{\dd^2y}{\dd t^2} = c_0(t)y.
\end{equation}
Then we can rewrite \eqref{eq:converted} in the rescaled time $s\coloneqq \varepsilon^{1/3}(\tau -t)$ so that it reads as a Airy equation of the form
\begin{equation}\label{eq:airy}
        \frac{\dd^2y}{\dd s^2} = sy,
\end{equation}
whose general solution is a linear combination of Airy functions of the first $\Ai(s)$ and second kind $\Bi(s)$ \cite[\S~$10.4$]{Abramowitz1965}.
By backward substitution, and in the past limit $t\to-\infty$, the pullback attractor solving \eqref{eq:riccati} is given by
\begin{equation}\label{eq:pullback_attractor}
        x_{\mathrm{pb}}(s) = -\varepsilon^{1/3}\frac{1}{\Ai(s)}\frac{\dd\Ai}{\dd s} .
\end{equation}
From \eqref{eq:pullback_attractor}, $x_{\mathrm{pb}}\to-\infty$ as $\Ai\to0$, which occurs at $s = s_0\approx 2.338$.
It thus follows that the tipping time of the pullback attractor is given by $\tau_{\mathrm{pb }} = \tau + s_0\,\varepsilon^{-1/3}$ and the opportunity window scales with the reciprocal of the cubic root of the parameter speed, i.e. $\Delta\tau\propto\varepsilon^{-1/3}$.

\subsection{Tipping time distribution}\label{subsec:exit}
We now consider the more general case of stochastic fluctuations of slow passage through the bifurcation by restoring $\sigma>0$ in \eqref{eq:model}.
As a result, the tipping time is treated as a random variable $\tau_{\mathrm{tip}}$ with an associated probability distribution. 
To characterize $\tau_{\mathrm{tip}}$ for \eqref{eq:model}, we follow the procedures outlined in \cite{Miller2012,Herbert2017}.
To start, consider again the singular limit $\varepsilon\to0$ of \eqref{eq:model}, written as a gradient system 
\begin{equation}\label{eq:gradient}
        \dd X_t=-V'\bigl(X_t,\mu\bigr)\,\dd t+\sqrt{2D}\,\dd W_t,
        \qquad
        V(x,\mu):=\mu x + \frac{x^3}{3},    
\end{equation}
where $V(x,\mu)$ is the potential of the vector field $f(x,\mu)$ and $D\coloneqq\sigma^2/2$.
The equilibria $x_{\pm}(\mu)=\pm\sqrt{-\mu}$ of $f$ correspond to the stationary points of $V$, respectively its local minimum and local maximum.
At any $\mu<0$, the potential barrier, separating the local minimum $x_+(\mu)$ from the local maximum $x_-(\mu)$ is
\begin{equation}\label{eq:barrier}
        \Delta V\coloneqq V\bigl(x_-(\mu)\bigr)-V\bigl(x_+(\mu)\bigr) = \frac43(-\mu)^{3/2} > 0.
\end{equation}
The critical transition of a solution of \eqref{eq:model}, starting in a neighbourhood of $x_+(\mu)$, can thus be associated to the escape of an overdamped particle across the potential barrier $\Delta V$ driven by noisy fluctuations of intensity $D$.
For fixed values of $\mu<0$ such that $\Delta V\gg D$, these escapes are rare events described by Freidlin--Wentzell theory \cite{Freidlin2012} which has been widely applied in many disciplines (see \cite{Santolin2025} and references therein).
In particular, the event of a trajectory $X_t$ of \eqref{eq:gradient}, starting near $x_+(\mu)$ and reaching $x_-(\mu)$ across the potential barrier $\Delta V$, satisfies a so-called \emph{large deviation principle}.
The \emph{first exit time} of such an event \cite[\S~$5.5$]{Gardiner2009}  
\begin{equation}\label{eq:first_exit}
        \tau_{\mathrm{exit}} \coloneqq \inf\{t\ge0\;:\;X_t\le x_-(\mu)\}, 
\end{equation}
is asymptotically exponentially distributed \cite[Thm.~$2.3$]{Berglund2013}
\begin{equation}\label{eq:large_deviation_principle}
        \tau_{\mathrm{exit}}\overset{d}{\to}\exp\bigl(k(\mu)\bigr) , \quad D\downarrow0,
\end{equation}
where $\overset{d}{\to}$ denotes convergence in distribution of a random variable and the rate parameter $k(\mu)$ is given by Kramers formula
\begin{equation}\label{eq:rate}
        k(\mu)=\frac{\sqrt{V''\bigl(x_+(\mu)\bigr)\,|V''\bigl(x_-(\mu)\bigr)|}}{2\pi}\,\exp\biggl(-\frac{\Delta V}{D}\biggr)
        =\frac{\sqrt{-\mu}}{\pi}\exp\biggl(-\frac{4(-\mu)^{3/2}}{3D}\biggr).
\end{equation}
\begin{figure}[t]
    \centering 
    \includegraphics[keepaspectratio, width=\textwidth]{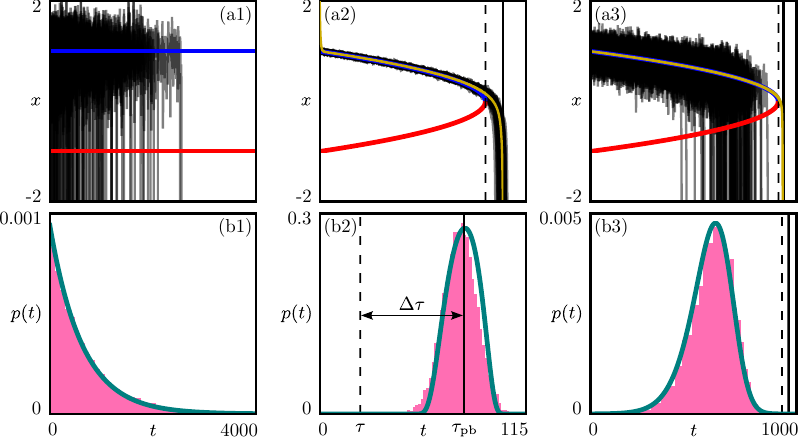}
    \caption{
      Comparison of the distributions of tipping times of the dynamic saddle--node normal form across different regimes: (a1)--(b1) noise-induced escapes at stationary parameter value $\mu=-1$ and noise level $D=0.25$; (a2)--(b2) slow--passage escapes at initial parameter value $\mu_0=-1$, speed $\varepsilon=10^{-2}$ and noise $D=10^{-3}$; (a3)--(b3) early escapes at initial parameter value $\mu_0=-1$, speed $\varepsilon=10^{-3}$ and noise $D=0.06$.
      In panels (a1)--(a3) we show Monte Carlo simulations of an ensemble of $10{,}000$ solutions of \eqref{eq:model} in each regime, reppresented by the black curves. 
      The blue and red curves reppresent the stable and unstable branches of equilibria $x_{\pm}(\mu(t))$ of $f$ respectively, while the yellow curve depicts the pullback attractor $x_{\mathrm{pb}}(t)$. 
      The vertical dahsed and solid black lines in (a2)--(a3) denote the breaking time $\tau$ and tipping time of the pullback attractor $\tau_{\mathrm{pb}}$ respectively. 
      In panels (b1)--(b3) the empirical distributions of the tipping time $\tau_{\mathrm{tip}}$ of the Monte Carlo simulations, depicted by the pink histograms, are compared with the analytical approximations \eqref{eq:large_deviation_principle}, \eqref{eq:small_noise} and \eqref{eq:large_noise} respectively, which are depicted by the teal curves.}
    \label{fig:tipping_times_distribution}
\end{figure}
Notice that the reciprocal of the escape rate in \eqref{eq:rate} gives the \emph{mean first exit time} of \eqref{eq:first_exit}, i.e. $\E[\tau_{\mathrm{exit}}] = 1/k(\mu)$; a comparison between the analytical distribution \eqref{eq:large_deviation_principle} and a finite sample Monte Carlo simulation is provided in Figure~\ref{fig:tipping_times_distribution}(b1).
This result, originally formulated for stationary diffusion \cite{Day1983} and jump processes \cite{Imkeller2006}, was then extended to the dynamic saddle-node normal form \cite{Miller2012}.
Specifically, in the limiting regime of low noise intensity $D\ll\varepsilon/2$, the approximate distribution of the tipping time $\tau_{\mathrm{tip}}$ has been shown to be \cite[Eq.~$(19)$]{Miller2012}
\begin{equation}\label{eq:small_noise}
        p_1(t)\coloneqq\frac{\varepsilon^{2/3}}{\sqrt{4\pi^{5}D}\,\Bi^{2}\bigl(-\mu(t)\varepsilon^{-2/3}\bigr)\,\mathcal J(t)^{1/2}}
        \exp\Biggl[-\frac{\varepsilon^{2/3}\,\Ai^{2}\bigl(-\mu(t)\varepsilon^{-2/3}\bigr)}
        {4\pi^{2}D\,\Bi^{2}\bigl(-\mu(t)\varepsilon^{-2/3}\bigr)\,\mathcal J(t)}\Biggr] ,
\end{equation}
with $\mathcal J(t)\coloneqq\int_{-\infty}^{t}\Ai^{4}\bigl(-\mu(s)\varepsilon^{-2/3}\bigr)\,\dd s$. 
We remark that $p_1(t)$ has support in the interval $\tau+ s_{\mathrm{inf}}\varepsilon^{-1/3} < t < \tau +  s_{\mathrm{sup}}\varepsilon^{-1/3}$, with $s_{\mathrm{inf}}\approx1.17$ and $s_{\mathrm{sup}}\approx3.27$, and that it contains the tipping time of the pullback attractor $\tau_{\mathrm{pb}} = \tau+s_0\,\varepsilon^{-1/3}$ since $s_{\mathrm{inf}} < s_0 < s_{\mathrm{sup}}$.
Furthermore \eqref{eq:small_noise} is approximately Gaussian locally to its maximum
\begin{equation}\label{eq:small_noise_approximation}
        \tau_{\mathrm{tip}}\overset{d}{\to}\mathcal{N}\biggl(\tau_{\mathrm{pb}}\bigl[1+O(D/\varepsilon)\bigr],O\bigl(D\varepsilon^{-5/3}\bigr)\biggr),
\end{equation}
with mean close to the tipping time of the deterministic solution, i.e. $\E[\tau_{\mathrm{tip}}]\approx\tau_{\mathrm{pb}}$ (see Figure~\ref{fig:tipping_times_distribution}(b2) for reference).

In the opposite regime of large noise intensity $\varepsilon/2\ll D$, and up to $-\mu(t)\gg D^{2/3}$, it was later shown that the approximate tipping time distribution has a maximum before the breaking time $\tau$ and it has the form \cite[Eq.~$(4)$]{Herbert2017}
\begin{equation}\label{eq:large_noise}
        p_2(t)\coloneqq k\bigl(\mu(t)\bigr)\,
        \exp\biggl[-\frac{D}{2\pi\varepsilon}\,
        \exp\biggl(-\frac{\Delta V\bigl(\mu(t)\bigr)}{D}\biggr)\biggr] ,
\end{equation}
where $k(\mu(t))$ is the Kramers escape rate given in \eqref{eq:rate} and $\Delta V(\mu(t))$ is the potential barrier given in \eqref{eq:barrier} (see Figure~\ref{fig:tipping_times_distribution}(b3) for reference).

The analytical approximations of the three different regimes prescribed in \cite{Berglund2013, Miller2012, Herbert2017} match the empirical observations obtained via Monte Carlo simulations as reported in Figure~\ref{fig:tipping_times_distribution}.
Notice how, in the small noise regime $D\ll \varepsilon/2$ (panel (b2)), almost all the mass of the tipping time distribution \eqref{eq:small_noise} is located after the breaking time $\tau$ and its mean approximates the tipping time of the pullback attractor $\tau_{\mathrm{pb}}$.
This setting thus prescribes the optimal scenario for which the (accurate) forecast of the breaking time $\tau$ provides valuable predictions for the estimation of the tipping time since it provides policymakers a window of opportunity $\Delta\tau$ to reverse the forcing of the system towards the saddle--node catastrophe.

\section{Forecasting problem}\label{sec:formulation}
From \S~\ref{sec:tipping}, we summarize three important results that motivate the forecasting of tipping times:
\begin{enumerate}
        \item for a slow passage through a saddle-node bifurcation the leading eigenvalue $\lambda(t)$ explicitly depends on the breaking time $\tau$ via Eq.~\eqref{eq:eigenvalue_time};
        \item the time at which a slowly ramped solution, in the noiseless limit $\sigma\to0$, tips $\tau_{\mathrm{pb}}$ is delayed with respect to the breaking time $\tau$ by a factor of $\varepsilon^{-1/3}$ that depends on the speed of the forcing, thereby creating an opportunity window $\Delta\tau$;
        \item in the small noise regime $\sigma=\sqrt{2D}\ll\sqrt{\varepsilon}$, the delayed tipping time of a noisy solution $\tau_{\mathrm{tip}}$ occurs almost always after the breaking time $\tau$, as per the support of the density $p_1(t)$ in Eq.~\eqref{eq:small_noise}, and on average close to the deterministic tipping time, as per Eq.~\eqref{eq:small_noise_approximation}.
\end{enumerate}
Suppose thus we are provided with a timeseries that is slowly ramped towards a saddle-node bifurcation, with noise $\sigma$ and speed $\varepsilon$ satisftying $0<\sigma^2\ll\varepsilon\ll1$. 
If we can accurately forecast the breaking time $\tau$ (i), we can then identify a window of opportunity (ii) and avoid the bifurcation-induced catastrophe with high probability (iii) by stopping or reversing the forcing of the parameter.
To formulate the forecasting problem of the breaking time $\tau$ we introduce the \emph{(squared) return rate} $\alpha(t)\coloneqq\lambda^2(t)$ which, from \eqref{eq:eigenvalue_time}, reads
\begin{equation}\label{eq:return_rate}
        \alpha(t) = 4\varepsilon(\tau-t)>0\quad\forall t<\tau . 
\end{equation}
Notice that, while $\alpha(t)\to0$ as $t\to\tau$, in the same fashion of $\lambda(t)\to0$ as $t\to\tau$, the return rate in \eqref{eq:return_rate} establishes an affine dependence of the return rate on time $t$ which benefits the forecasting of $\tau$ using timeseries data.

To understand this suppose that $\mathcal{X}_1\coloneqq\{X^{(1)}_j\coloneqq X_{t^{(1)}_j}\}_{j=1,\dots,n}$, $t^{(1)}_j<t^{(1)}_{j+1}$, is a solution of \eqref{eq:model} that is at safe distance from the bifurcation, i.e. $t^{(1)}_n\ll\tau$.
Assume also that we can estimate the true return rate \eqref{eq:return_rate} from the sample $\mathcal{X}$, which we denote with $\hat\alpha_1\coloneqq\hat\alpha(t^{(1)}_n)$.
We are then presented with another sample $\mathcal{X}_2$, which is also at safe distance from the bifurcation but closer to it than the previous sample $\mathcal{X}_1$, i.e. $t^{(1)}_n<t^{(2)}_n$; once again, assuming we can do so, we estimate the return rate $\hat\alpha_2$.
Suppose we repeat this procedure for $N$ different samples $\mathcal{X}_1,\dots,\mathcal{X}_N$, yielding a sequence $\hat{\alpha}\coloneqq\{\hat\alpha_j\}_{j=1,\dots,N}$.
Given this sequence $\hat\alpha$ we can then fit a linear model of the form $\hat\alpha(t) = \hat m t + \hat q$, with slope $\hat m$ and intercept $\hat q$. 
The forecasted breaking time $\hat\tau$ thus corresponds to the zero-crossing of the linear fit, which by extrapolation reads $\hat\tau = -\hat q/\hat m$.

Consequentially, we translate the problem of forecasting the breaking time $\tau$ from timeseries data $\mathcal{X}_1,\dots,\mathcal{X}_N$ into a simple, linear extrapolation $\hat\tau$ from estimated return rates $\hat\alpha$.
This procedure, however, features two bottlenecks: first, we need $\hat\alpha$ to be a sequence of precise and accurate estimators of the ground truth $\alpha$ and; second, given those estimates $\hat\alpha$, we need the linear extrapolation $\hat\tau=-\hat q/\hat m$ to be as close as possible to the ground truth $\tau$ to not shrink the window of opportunity.
In particular, depending on applications, we would prefer to slightly underestimate the breaking time $\hat\tau\lesssim\tau$, rather than slightly overestimating it $\hat\tau\gtrsim\tau$, thereby biasing our forecast to minimize missed alarms (i.e. false negatives) in favour of spurious warnings (i.e. false positives).
In the remainder of this section we present three data-driven techniques to accurately extract the estimates $\hat\alpha$ from a sample $\mathcal{X}$ while we address the extrapolation problem in \S~\ref{sec:observation}.

\subsection{Ornstein-Uhlenbeck approximation}\label{subsec:OUP}
Let us start by characterizing the statistical properties of the noisy fluctuations of a solution $\mathcal{X}$ of \eqref{eq:model} around the pullback attractor $x_{\mathrm{pb}}(t)$ given in Eq.~\eqref{eq:pullback_attractor}.
We write the solution variable as a superposition $X_t\coloneqq x_{\mathrm{pb}}(t)+Z_t$ of a deterministic $x_{\mathrm{pb}}(t)$ and a stochastic component $Z_t$.
Substituting this into \eqref{eq:model}, recalling that by definition $\dd x_{\mathrm{pb}}=f(x_{\mathrm{pb}},\mu)\dd t$, we get
\begin{equation}\label{eq:fluctuations}
  \dd Z_t=\bigl(-2 x_{\mathrm{pb}}(t)\,Z_t-Z_t^2\bigr)\dd t+\sigma\,\dd W_t .
\end{equation}
At any time $t\ll\tau$ away from the saddle-node bifurcation $f(x_{\mathrm{pb}},\mu)\approx f\bigl(x_+(\mu),\mu\bigr)$ which implies $\partial_x f(x_{\mathrm{pb}},\mu)\approx\lambda(t)$.
Substituting this into \eqref{eq:fluctuations} gives
\begin{equation}\label{eq:fluctuations_scaled}
        \dd Z_t=\bigl(\lambda(t)\,Z_t-Z_t^2\bigr)\dd t+\sigma\,\dd W_t .
\end{equation}
It can be shown that $|Z_t|\approx\sigma/\sqrt{2|\lambda(t)|}$ \cite[Thm.~$4.1$]{Kuehn2011}, and thus $Z_t^2=o\bigl(|\lambda(t)Z_t|\bigr)$ for any $\sigma\ll\sqrt{\varepsilon}$ and any $t\ll\tau$ such that $\sigma\ll\sqrt{2}\bigl(4\varepsilon(\tau-t)\bigr)^{3/4}$.
We can therefore reasonably approximate \eqref{eq:fluctuations_scaled} by neglecting the nonlinear term $Z_t^2$, thereby obtaining 
\begin{equation}\label{eq:oup}
        \dd Z_t\approx\lambda(t)\,Z_t\dd t+\sigma\,\dd W_t .
\end{equation}
The linear, nonautonomous Eq.~\eqref{eq:oup} is a time-inhomogeneous Ornstein-Uhlenbeck process.
Suppose the sample $\mathcal{X}$ is recorded across a time window $T_w>0$ in which $\lambda(t)\approx\lambda$ is approximately constant. 
Then such sample can be considered as an approximate solution of the stationary process \eqref{eq:oup} which admits a Gaussian distribution \cite[\S~$3.8$]{Gardiner2009} with zero mean $\E[Z_t]\eqqcolon\bar{Z}=0$ and variance $\Var[Z_t]\eqqcolon\gamma=\sigma^2/2|\lambda|$.
Written in terms of the return rate $\alpha=\lambda^2$, the last identity reads
\begin{equation}\label{eq:oup_variance}
        \gamma = \frac{\sigma^2}{2\sqrt{\alpha}} .
\end{equation}
Letting $\hat\gamma$ be the sample variance of $\mathcal{X}$, then \eqref{eq:oup_variance} gives a sample estimate of the return rate
\begin{equation}\label{eq:OUP}
        \hat\alpha = \frac{\sigma^4}{4\hat\gamma^2} .
\end{equation}
The use of the sample variance as an early-warning signal dates back to the very first attempts in predicting bifurcation-induced tipping points and it has since been established as a benchmark precursor of critical slowing down (see e.g. \cite{Kuehn2011} and references therein).

\subsection{Autoregressive approximation}\label{subsec:AC1}
In the same time window $T_w>0$ above, the homogeneous Ornstein-Uhlenbeck process has correlation function $\Corr(Y_t,Y_{t+Tw})=\gamma\exp\bigl(-\sqrt{\alpha}T_w\bigr)$ \cite[\S~$3.8$]{Gardiner2009}, with $\gamma$ given in \eqref{eq:oup_variance}.
As a result, the discrete observations in the sample $\mathcal{X}$ can be approximated by fitting a continuous-time autoregressive model of order one
\begin{equation}\label{eq:autoregressive}
        X_{j+1}\approx\rho X_{j} + \sqrt{\gamma(1-\rho^2)}\,\eta_t , \quad \eta_t\sim\mathcal{N}(0,1),
\end{equation}
where 
\begin{equation}\label{eq:ac1}
        \rho\coloneqq\gamma^{-1}\Corr(X_t,X_{t+\Delta}), \quad \Delta=t_{j+1}-t_j, 
\end{equation}
is the autocorrelation coefficient at lag$-1$.
Given the sample mean $\hat X$ the sample autocorrelation coefficient can be obtained by maximum likelihood estimation \cite{Sahalia2002} and it reads
\begin{equation}\label{eq:ac1_sample}
        \hat\rho = \frac{\sum_{j=1}^{n-1}(X_j - \hat X)(X_{j+1} - \hat X)}{\sum_{j=1}^{n}(X_j - \hat X)^2}.
\end{equation}
Using \eqref{eq:ac1_sample} one can estimate the return rate of the sample $\mathcal{X}$ using the autoregressive approximation 
\begin{equation}\label{eq:AC1}
        \hat\alpha = \biggl(\frac{\log\hat\rho}{\Delta}\biggr)^2.
\end{equation}
Autoregressive based early-warning signals have found widespread applications in ecology \cite{Dakos2012} and climate science \cite{Boers2021}.
In particular, estimators of the form \eqref{eq:AC1} have been used in \cite{Ditlevsen2023} for estimating the tipping time of the subpolar gyre of the Atlantic Meridional Overturning Circulation system, and in \cite{Ashwin2025} for establishing a framework that quantifies the predicting skill of early-warning signals.

\subsection{Euler-Maruyama approximation}\label{subsec:LLS}
A solution of the nonlinear system \eqref{eq:model} in a time window $T_w>0$ as per above, can be approximated by the process \cite[\S~$8.4$]{Lord2014}
\begin{equation}\label{eq:EM}
        X_{t+1} = X_t + f(X_t,\mu)\Delta + \sigma\sqrt{\Delta}\,\xi_t , \quad \xi_t\sim\mathcal{N}(0,1) .
\end{equation}
The sequence above is the Euler-Maruyama approximation of \eqref{eq:model} and it corresponds to the stochastic generalization of the Euler time-stepping method for ordinary differential equations.
We can rearrange the terms in \eqref{eq:EM} to define the first-order time increments of the Euler-Maruyama approximation
\begin{equation}\label{eq:increments}
        Y_t\coloneqq\frac{X_{t+1} - X_t}{\Delta} = f(X_t,\mu) + \tilde{\xi}_t , \quad \tilde{\xi}_t\sim\mathcal{N}\biggl(0,\frac{\sigma^2}{\Delta}\,\biggr) .
\end{equation}
Provided a sample $\mathcal{X} = \{X_j\}_{j=1,\dots,n}$ in a time window $T_w=n\Delta$, and a choice of the model function $f$, Eq.~\eqref{eq:increments} can be interpreted as an ordinary least--squares regression problem, with normally distributed residuals $r_j\coloneqq Y_{j+1} - f(X_j)=\tilde{\xi}_t$.
This is part of a broader class of drift inference methods called \emph{quasi-maximum likelihood estimators} \cite{Jimenez2020}.

Away from the saddle-node bifurcation, the sample is approximately an Ornstein-Uhlenbeck process of the form \eqref{eq:oup} and thus we can reasonably fit a linear model of the form $f(x;\theta) = \theta x$ to the time increments $Y_j$.
Substituting this ansatz into \eqref{eq:increments}, and minimizing the square root of the sum of the squared residuals across the sample, yields the optimization problem
\begin{equation}\label{eq:ordinary_least_squares}
        \hat\theta = \argmin_{\theta\in\mathbb{R}}\Biggl(\sqrt{\sum_{j=1}^{n-1}\bigl(Y_j-f(X_j;\theta)\bigr)^2}\Biggr) ,
\end{equation}
which is linear in $\theta$.
The solution of \eqref{eq:ordinary_least_squares} is given by the normal equation
\begin{equation}\label{eq:normal_equation}
        \hat\theta = (A^TA)^{-1}A^TY,
\end{equation}
where $Y=(Y_j)_{j=1,\dots,n-1}$ is the vector of observed time increments of the sample $\mathcal{X}$ and $A=(X_j)_{j=1,\dots,n-1}$ is the model vector of the linear fit $f$ to the sample $\mathcal{X}$.
Notice that, by definition, $\hat{\theta}$ is an estimate of the leading eigenvalue $\lambda$ of the underlying nonlinear evolution and thus
\begin{equation}\label{eq:LLS}
        \hat\alpha = \hat{\theta}^{\,2},
\end{equation}
is a sample estimate of the return rate using an ordinary (linear) least--squares approximation. 
The application of quasi-maximum likelihood estimations for the construction of early-warning signals is much more recent than the more established autoregressive and Ornstein-Uhlenbeck approximation. 
We recently implemented this approach in \cite{Papapicco2026} to construct an interpretable early--warning indicator adapted from Kramers escape formula \eqref{eq:rate}.

\section{Observation protocol}\label{sec:observation}
The return rate estimators using the Ornstein--Uhlenbeck \eqref{eq:OUP}, autoregressive \eqref{eq:AC1} and Euler--Maruyama \eqref{eq:LLS} approximations provide three different, data--driven pipelines to extract $\hat\alpha$ from a given sample $\mathcal{X}$ collected over a time window $T_w=n\Delta$, with $n>0$ number of observations sampled at (regular) frequency $\Delta$.
As noted at the beginning of \S~\ref{sec:formulation}, however, we then use these estimates to forecast the breaking time $\hat\tau$ given a sequence of samples $\mathcal{X}_1,\dots,\mathcal{X}_N$.
In this section we address this problem in two parts.
First, we discuss the accuracy and precision associated to the estimators of the return rate $\hat\alpha$ proposed above.
Second, we formulate the extrapolation problem of $\hat\tau$ given a sequence $\hat\alpha=\{\hat\alpha_j\}_{j=1,\dots,N}$ of estimators of the return rate and construct the full, data-driven pipeline for the accurate forecasting problem given input timeseries data.

\subsection{Estimating the return rate}\label{subsec:estimators}
The optimal estimation of the return rate from a finite sample is, perhaps unsurprisingly, highly limited by the quality of the input data provided by the samples $\mathcal{X}_1,\dots,\mathcal{X}_N$.
We hereby discuss three major design constraints that modellers and scientists have to consider when applying the proposed forecasting methodology to their samples.

\subsubsection{Admissible time horizon}\label{subsubsec:horizon}

First observe that the analytical return rate \eqref{eq:return_rate}, on which the sample estimators \eqref{eq:OUP}, \eqref{eq:AC1} and \eqref{eq:LLS} are based on, is obtained from the time--dependent eigenvalue \eqref{eq:eigenvalue_time}, which is itself evaluated along the drifting quasi--steady equilibrium $x_+(t)$.
However, if the sample $\mathcal{X}$ is assumed to be a slowly forced solution of \eqref{eq:model} away from the bifurcation, then its trend follows the pullback attractor \eqref{eq:pb_attractor} instead of $x_+(t)$.
In the slow regime (after the fast relaxation of the initial condition) the distance at which $x_{\mathrm{pb}}(t)$ tracks $x_+(\mu(t))$ increases with time and it reaches its maximum at the saddle-node bifurcation $|x_{\mathrm{pb}}(t)-x_+(\mu(t))|=O(\varepsilon^{1/3})$, $t\to\tau$ \cite[Thm.~$2.2.2$]{Berglund2006}.
We are thus interested in quantifying a time horizon $t=T_c$ for which the return rate $\alpha(t)$ of the quasi-steady equilibrium is (relatively) accurate to that of the pullback attractor $\alpha_{\mathrm{pb}}(t)=\lambda_{\mathrm{pb}}^2(t)=4x_{\mathrm{pb}}^{\,2}(t)$.

Sufficiently far from the saddle-node bifurcation we can exploit the asymptotic expansion of Airy's function \cite[\S~$10.4$]{Abramowitz1965} to approximate \eqref{eq:pullback_attractor} to first-order as
\begin{equation}\label{eq:asymptotic_pullback_attractor}
  x_{\mathrm{pb}}(s)\approx \varepsilon^{1/3}\biggl(\sqrt{s} + \frac{1}{4s}\biggr) , \quad s=\varepsilon^{1/3}(\tau-t)\to\infty.
\end{equation}
In this regime, the return rate for the pullback attractor is thus
\begin{equation}\label{eq:pb_return_rate}
  \alpha_{\mathrm{pb}}(t) = \alpha(t) + \frac{1}{4(\tau-t)^2} + 2\sqrt{\frac{\varepsilon}{\tau-t}},
\end{equation}
and dividing each term in \eqref{eq:pb_return_rate} by $\alpha(t)$ gives
\begin{equation}\label{eq:lag_correction}
  \frac{\alpha_{\mathrm{pb}(t)}}{\alpha(t)} = 1 + \frac{1}{16 \varepsilon(\tau - t)^3} + \frac{1}{2\sqrt{\varepsilon(\tau - t)^3}} = 1 + 2\eta(t) + \eta^2(t) = \bigl(1+\eta(t)\bigr)^2 ,
\end{equation}
where $\eta(t)\coloneqq 1/(4\sqrt{\varepsilon(\tau-t)^3})$ quantifies the lag correction of the return rate of the pullback attractor.

\begin{figure}[t]
    \centering 
    \includegraphics[keepaspectratio, width=\textwidth]{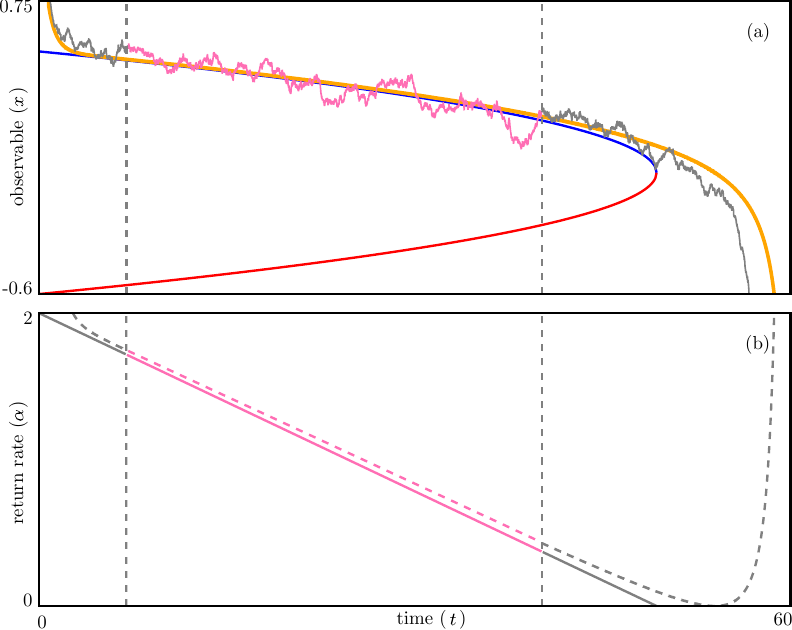}
    \caption{Illustration of the correction in the return rate $\alpha(t)$ between the quasi--steady equilibrium $x_+(t)$ and the pullback attractor $x_{\mathrm{pb}}(t)$.
      In panel (a) we show a solution of \eqref{eq:model}, given by the grey timeseries, tracking the pullback attractor, depicted by the yellow curve, and eventually separating from the drifting quasi--steady equilibria, given by the blue curve. 
      This solution is generated by setting $\varepsilon=10^{-2}$ and $\sigma=0.6\sqrt{\varepsilon}$ in \eqref{eq:model}.
    In panel (b) we show the effect of this separation on the return rate $\alpha$; the dashed curve is the return rate evaluate along the pullback attractor, while the solid curve is the return rate evaluated at the quasi--steady equilibrium.
    The vertical dashed grey lines in both panels represent the admissible time horizon for which a subseries of the solution can be used for forecasting, highlighted in pink.
  }
    \label{fig:admissible_time_horizon}
\end{figure}

From \eqref{eq:lag_correction} it follows that $\alpha(t)\approx \alpha_{\mathrm{pb}}(t)$ at any time $t\ll\tau$ such that $\eta(t)\ll1$, from which we obtain the admissible time horizon
\begin{equation}\label{eq:admissible_time_horizon}
  T_c \lesssim\tau - \varepsilon^{-1/3} .
\end{equation}

The critical threshold \eqref{eq:admissible_time_horizon} gives us a criterion for which the estimation of the frozen return rate yields admissible results for the extrapolation of the tipping time, as the lag correction $\eta(t)$ of the pullback attractor with respect to the quasi-steady equilibrium $x_+(t)$ is negligible. 
This is depicted in Figure~\ref{fig:admissible_time_horizon} in which a sample path, shown by the grey timeseries in panel (a), follows the pullback attractor solution $x_{\mathrm{pb}}(t)$ of \eqref{eq:model}, shown by the yellow curve. 
Notice how $x_{\mathrm{pb}}(t)$ tracks the drifting quasi--steady equilibrium $x_+(t)$, shown by the blue curve, for some time before the separation becomes visible close to the saddle--node bifurcation.
This is reflected in the difference in the return rates of $x_{\mathrm{pb}}(t)$ and $x_+(t)$, depicted by the dashed and solid grey curves in panel (b) respectively.
This discrepancy thus identifies an admissible time horizon, shown by the vertical dashed grey curves, bounded below by the exctintion of the transient of the pullback attractor, converging to the (attracting) slow manifold of the nonautonomous system, and above by \eqref{eq:admissible_time_horizon}.

\subsubsection{Biased estimators}\label{subsubsec:effects}
Provided that the sample is in the admissible time horizon $t<T_c$ quantified above, a local estimation of $\alpha(t)$ requires that the sample $\mathcal{X}$ is recorded over a time window $T_w$ in which the leading eigenvalue $\lambda(t)$ at the quasi--steady equilibrium remains approximately constant.
This in turn requires the sample to be limited in a finite observation window $T_w=n\Delta$ that induces biases and uncertainties in the estimation of the return rate.

To quantify those, we shall henceforth adopt an ensemble approach in which $\Omega = \{\mathcal{X}_{j}\}_{j=1,\dots,N}$ denotes a collection of $N>1$ independent and identically distributed finite samples that correspond to different noise realizations of \eqref{eq:model}.
In particular, let us assume that each sample $\mathcal{X}_j$ in the ensemble $\Omega$ is a sequence of stationary residuals obtained by detrending a finite windowed solution of \eqref{eq:model} so that $t^{(1)}_n\equiv t^{(N)}_n$.
Then, the resulting Ornstein-Uhlenbeck processes \eqref{eq:oup} can be approximated by order one autoregressive models of the form \eqref{eq:autoregressive}.

Then let $\hat \beta$ be a generic random variable estimating some ground truth $\beta$ from a single sample in the ensemble.
Denote $\E[\hat \beta]$ to be the ensemble mean of $\hat\beta$ over $\Omega$ and $\Var[\hat\beta]$ to be the variance of $\hat\beta$ over $\Omega$.
We thus define the bias $b_{\hat\beta}$ and uncertainty $u_{\hat\beta}$ of the estimator $\hat\beta$ in the canonical fashion
\begin{equation}\label{eq:bias_and_uncertainty}
  b_{\hat\beta}\coloneqq\E[\hat\beta - \beta] \quad \text{and}\quad u_{\hat\beta}\coloneqq\sqrt{\Var[\hat\beta]} ,
\end{equation}
where, depending on the context, these are also referred to as accuracy and precision respectively; here these two terminologies will be used interchangeably.

\paragraph{Bias in the sample variance}
The return rate estimator \eqref{eq:OUP} uses the sample variance in the Ornstein-Uhlenbeck approximation.
It is a standard result that the sample variance $\hat\gamma=\Var[\mathcal{X}]$ is a biased estimator for autocorrelated data.
Specifically, its accuracy and precision at lag-$1$ are given by \cite[Eq.~$(4.6)$]{Law1991}
\begin{equation}\label{eq:bias_law}
  b_{\hat\gamma} = -\frac{2\gamma\rho}{(n-1)(1-\rho)}\quad \text{and} \quad u_{\hat\gamma} = \sqrt{\frac{2\gamma^2(1+\rho^2)}{(n-1)(1-\rho^2)}},
\end{equation}
where $\rho$ is the autocorrelation coefficient given in \eqref{eq:ac1}.
Following a standard delta method \cite[\S~$5.5.4$]{Casella2024}, and using \eqref{eq:oup_variance}, we obtain a second-order approximation of the bias in the estimated leading eigenvalue, which reads
\begin{equation}\label{eq:taylor_expansion}
  b_{\hat\lambda}^{\text{\tiny OUP}} = \frac{|\lambda|}{\gamma}b_{\hat\gamma} - \frac{|\lambda|}{\gamma^2}u^2_{\hat\gamma}.
\end{equation}
Substituting \eqref{eq:bias_law} into \eqref{eq:taylor_expansion} yields a closed form expression for the bias in the estimation of the leading eigenvalue of a finite sample of residuals given its sample variance
\begin{equation}\label{eq:OUP_bias}
  b_{\hat\lambda}^{\text{\tiny OUP}} = -\frac{2|\lambda|(1+\rho+2\rho^2)}{(n-1)(1-\rho^2)}.
\end{equation}

\paragraph{Bias in the sample autocorrelation}
The return rate estimator \eqref{eq:AC1} uses the squared logarithm of the sample autocorrelation \eqref{eq:ac1_sample}, whose bias is quantified in \cite[Eq.~$(21)$]{Kendall1954}
\begin{equation*}
  b_{\hat\rho} = - \frac{1 + 3\rho}{n-1} .
\end{equation*}
Applying again the delta method, the above yields the bias in the leading eigenvalue \cite[Eq.~$(12)$]{Yu2012}
\begin{equation}\label{eq:AC1_bias}
  b_{\hat\lambda}^{\text{\tiny AC1}} = - \frac{3 + \rho^{-2}}{2\Delta(n-1)} .
\end{equation}
\begin{figure}[t]
    \centering 
    \includegraphics[keepaspectratio, width=\textwidth]{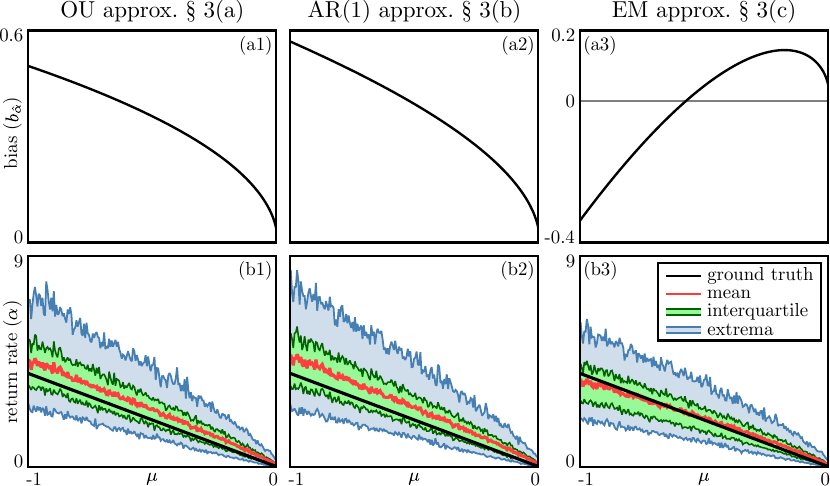}
    \caption{Numerical validation of the return rate bias \eqref{eq:return_rate_bias} given by: (a1)--(b1) the sample variance estimator \eqref{eq:OUP}; the autoregressive model fit \eqref{eq:AC1} and; the Euler--Maruyama approximation \eqref{eq:LLS}.
      In panels (a1)--(a3) we plot \eqref{eq:return_rate_bias}, in the parameter range $\mu\in[-1,0)$, using \eqref{eq:OUP_bias}, \eqref{eq:AC1_bias} and \eqref{eq:LLS_bias} respectively.
      In panels (b1)--(b3) we plot the results of Monte Carlo simulations of an ensemble $\Omega=\{\mathcal{X}_j\}_{j=1,\dots,N}$ of $N=200$ independent but stationary solutions of \eqref{eq:model} with $n=400$ realizations and noise level $\sigma=0.02$.
      At discrete parameter values $\mu\in[-1,0)$, and for each sample $\mathcal{X}_j\in\Omega$, we estimate $\hat\alpha$ using the three techniques discussed above; this forms an ensemble of estimators from which representative statistics are computed.
    In particular the ensemble mean, shown by the red curves, are compared to the ground truth, black curves, to assess the bias in the estimators; the ensemble interquartile range, shown by the green shaded regions, are used to quantify the uncertainty in the estimators while the blue shaded regions provide the $5^{\text{th}}$-$95^{\text{th}}$ range of the estimators.}
    \label{fig:bias_estimators}
\end{figure}

\paragraph{Bias in the ordinary least--squares solution}
The return rate estimator \eqref{eq:LLS} is the square of the ordinary least--squares solution of the Euler-Maruyama approximation \eqref{eq:normal_equation}, which coincides with the estimated leading eigenvalue of the sample $\mathcal{X}$ by definition.
Its bias is thus more straightforward to compute \cite[Eq.~$(4.5)$]{Florens1989} and it reads
\begin{equation}\label{eq:LLS_bias}
  b_{\hat\lambda}^{\text{\tiny OLS}} = \frac{\Delta\lambda^2}{2} - \frac{1+3\rho}{\Delta(n-1)}. 
\end{equation}
Notice that, while \eqref{eq:OUP_bias} and \eqref{eq:AC1_bias} are negative for all samples up to the saddle-node bifurcation $\lambda=-2\sqrt{-\mu}<0$, the bias in the least--squares estimator \eqref{eq:LLS_bias} changes sign at 
\begin{equation}\label{eq:unbiased_LLS}
  \mu = - \frac{1+3\rho}{2\Delta^2(n-1)} ,
\end{equation}
for which \eqref{eq:LLS} becomes the only unbiased estimator of the leading eigenvalue.

\paragraph{Accuracy and precision of the estimators}
Eqs.~\eqref{eq:OUP_bias}-\eqref{eq:LLS_bias} measure the accuracy for the estimation of the leading eigeinvalue $\lambda$ of a sample. To convert those into biases of the return rate $\alpha$, notice that
\begin{equation}
  \E[\hat\alpha] = \E[\hat\lambda^2] = \E[\hat\lambda]^2 + \Var[\hat\lambda] = (b_{\hat\lambda} + \lambda)^2 + u^2_{\hat\lambda} ,
\end{equation}
where in the last equality we used the definitions \eqref{eq:bias_and_uncertainty}.
Expanding the quadratic term and using the definition of the return rate $\alpha$ we obtain
\begin{equation}\label{eq:return_rate_bias}
  b_{\hat\alpha} = b^2_{\hat\lambda} + 2\lambda b_{\hat\lambda} + u^2_{\hat\lambda},
\end{equation}
from which we conclude that the accuracy of the estimators depend on both the accuracy $b_{\hat\lambda}$  and the precision $u_{\hat\lambda}$ of the leading eigenvalue $\lambda$.
In the high sampling frequency limit this uncertainty is approximated by the Cramer--Rao inequality \cite[Thm.~$7.3.9$]{Casella2024} 
\begin{equation}\label{eq:cramer_rao_bound}
  u_{\hat\lambda}^2 \to \frac{2|\lambda|}{\Delta(n-1)},\quad\Delta\to 0. 
\end{equation}

Monte Carlo simulations, depicted in Figure~\ref{fig:bias_estimators}(b1--b3), show that estimating the return rate of a sample $\mathcal{X}$ using the Euler--Maruyama approach \eqref{eq:LLS} yields a more accurate and more precise estimation than the established sample variance \eqref{eq:OUP} and autocorrelation \eqref{eq:AC1}.
This result matches the second order approximations for the bias of $\hat\alpha$ derived in \eqref{eq:OUP_bias}, \eqref{eq:AC1_bias} and \eqref{eq:LLS_bias}, which we depict across the parameter interval $\mu[-1,0)$ in Figure~\ref{fig:bias_estimators}(a1--a3).

\subsubsection{Optimal observation length}\label{subsubsec:window}
The topic of selecting an optimal window length to forecast critical transitions was addressed in \cite{Chen2022}.
In particular, the authors showed that windowed observations with less than $10\%$ of the data from a given sample yield spurious early-warning signals.
Such a study however only considers fixed window lengths across a solution approaching a critical transition, which neglects the changing characteristics of the solution in the slow passage through the saddle-node.
Here we derive the optimal, time-dependent window length of a sample for the accurate forecasting of the tipping time.

\begin{figure}[!t]
    \centering 
    \includegraphics[keepaspectratio, width=\textwidth]{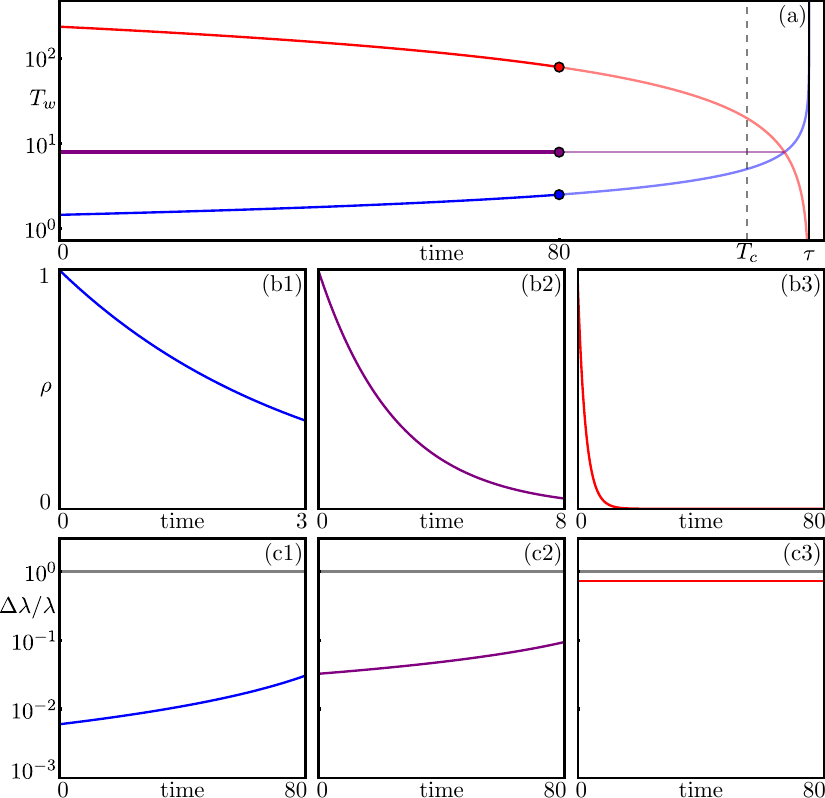}
    \caption{Visualization of the bounds of the observation length and its effect on the quality of the samples for the forecasting problem.
      Panel (a), vertical axis in logarithmic scale: the optimal length $T_w$ (purple line) for a dynamic saddle--node normal form \eqref{eq:model} is fixed over time to be within the lower \eqref{eq:lower_bound} and upper bound \eqref{eq:upper_bound} (blue and red curves respectively). The time horizon \eqref{eq:admissible_time_horizon} and breaking time \eqref{eq:breaking_time} are indicated for reference (dashed grey and solid black vertical line respectively).
      Panels (b1--b3): decay of the autocorrelation coefficient \eqref{eq:ac1} within the window length provided by the colored circles in panel (a). Notice how an observation length on the lower bound (panel (b1)) corresponds to shorter time of the relaxation time $1/|\lambda|$.
      The three panels have different time scales to emphasize the lethargic effect of relaxation of the autocorrelation coefficient $\rho$.
      Panels (c1--c3), vertical axes in logarithmic scale: relative change in the eigenvalue across the window length corresponding to the colored curves in panel (a). Notice how an observation length on the upper bound (panel c3) corresponds to near unity relative changes at all time, i.e. $\Delta\lambda\approx|\lambda|$.}
    \label{fig:observation_window}
\end{figure}

Given a sample $\mathcal{X}=\{X_{t_1},\dots,X_{t_n}\}$ recorded across a window $T_{w}=t_n-t_1$, we look for bounds of the form
\begin{equation}\label{eq:bounds}
  L(t)\ll T_w \ll U(t) ,
\end{equation}
in which the residuals retain the necessary properties for the estimation of the return rate \eqref{eq:return_rate}.
The lower bound $L(t)$ is derived by imposing that the sample $\mathcal{X}$ is observed on a window that is much longer than the relaxation time $1/|\lambda|$ which, from \eqref{eq:eigenvalue_time}, reads
\begin{equation}\label{eq:lower_bound}
  L(t) = \frac{1}{2\sqrt{\varepsilon(\tau-t)}}\ll T_w.
\end{equation}
This requirement ensures that the realizations $X_j$ in the sample are independent fluctuations around the stable equilibrium $x_+(t)$ due to the decay in autocorrelation. 
In Figure~\ref{fig:observation_window} we can see the effects of this choice; when the observation window is too short (i.e. $T_w=L(t)=3$) the autocorrelation coefficient does not decay enough to guarantee independence in the observed fluctuations around the drifting quasi--steady equilibrium (panel (b1)).
On the other hand, when the observation length is sufficiently large (i.e. $T_w = 8 > L(t)$) the same coefficient decays substantially more (at about $5\%$ of its starting value) thereby allowing enough time for the stochastic perturbations to relax back to $x_s(t)$ (panel (b2--b3)).

Regarding the upper bound $U(t)$, it is derived by imposing that the sample $\mathcal{X}$ is quasi--stationary.
To derive $U(t)$ observe that, from \eqref{eq:eigenvalue_time}, the instantaneous rate of change of the eigenvalue is
\begin{equation*}
  \frac{\dd \lambda}{\dd t} \eqqcolon\dot\lambda(t) = \frac{2\varepsilon}{|\lambda(t)|} ,
\end{equation*}
which is positive and monotonically increasing for any $t<\tau$, which implies that $\dot\lambda(t_n)=\sup_{t_1<t<t_n}\dot\lambda(t)$ over any finite time window $(t_1,t_n)$.
By setting $T_{w}=t_n-t_1$ it follows
\begin{equation*}
  0<\Delta\lambda = \lambda(t_n) - \lambda(t_1) = \int_{t_1}^{t_n}\dot\lambda(t)dt < \dot\lambda(t_n)T_{w}  = \frac{2 \varepsilon}{|\lambda(t_n)|}T_{w},
\end{equation*}
and dividing each term by $|\lambda(t_n)|$ yields
\begin{equation*}
  \frac{\Delta\lambda}{|\lambda(t_n)|} < \frac{2 \varepsilon}{\lambda^2(t_n)}T_{w} .
\end{equation*}
Our requirement for the sample $\mathcal{X}$ to be quasi-stationarity across $T_{w}$ corresponds to small changes of the leading eigenvalue compared to its endpoint value, i.e. $\Delta\lambda/|\lambda(t_n)|\ll 1$, which is guaranteed by imposing
\begin{equation*}
  \frac{2 \varepsilon}{\lambda^2(t_n)}T_{w} \ll 1 ,
\end{equation*}
yielding a sufficient condition for the upper bound of the window length
\begin{equation}\label{eq:upper_bound}
  T_w\ll U(t) = \frac{\lambda^2(t)}{2 \varepsilon} = 2(\tau - t),
\end{equation}
where the last equality follows from \eqref{eq:return_rate}.
Once again, in Figure~\ref{fig:observation_window} we depict how the upper bound constraint affects the stationarity of the sample. 
At a given $\varepsilon>0$, if the sample is too long it leads to changes in the leadeing eigenvalue that are comparable to the values of the same quantity at the end of the window, thereby not ensuring quasi--stationarity (panel (c3)).
For short enough observation lengths however, the same quantity reduces to a negligible fraction, which guarantees that the drift of $x_s(t)$ is not overly pronounced (panels (c1--c2)).

Is thus clear that the optimal observation length of a striding window must strike a balance between guaranteeing conditional independence of the collected obseravations while ensuring quasi--stationarity of the sample.
Substituting \eqref{eq:lower_bound} and \eqref{eq:upper_bound} into \eqref{eq:bounds} provides us such window length protocol
\begin{equation}\label{eq:window_length}
  \frac{1}{2\sqrt{\varepsilon(\tau-t)}}\ll T_w \ll 2(\tau - t),
\end{equation}
where notice that only the lower bound depends on the speed of the parameter $\varepsilon$.
We further observe that, since both $L(t)$, $U(t)$ depend on $t$ through its distance to the breaking time $\tau$, it may be necessary to vary the window length $T_w(t)$ as a function of time in order for it to be bounded within the optimality region. 
However, since $L(t)$ is monotonically increasing and $U(t)$ is monotonically decreasing in $t$, it is always possible to find a fixed value of $T_w(t)=T_w$ that can be used for all times up to the admissible time horizon $T_c$.
An example of this property is shown by the purple line in Figure~\ref{fig:observation_window}(a) bounded from above by the red line, depicting the decreasing $U(t)$, and from below by the blue curve, depicting the increasing $L(t)$, for all $t<T_c$.

\subsection{Extrapolating the tipping time}\label{subsec:extrapolators}
Having derived the design conditions for the collection of samples that yield optimal estimators $\hat\alpha$ we are now ready to specify the full data-driven pipeline for the forecasting of the breaking time $\tau$.
As anticipated in \S~\ref{sec:formulation}, given a sequence of estimators $\hat\alpha=\{\hat\alpha_j\}_{j=1,\dots,N}$, the forecasting of $\tau$ is translated into an extrapolation problem $\hat\tau = -\hat q/\hat m$ via the linear fit $\hat\alpha(t)=\hat m t + \hat q$ of slope $\hat m$ and intercept $\hat q$.
As for the Euler--Maruyama approximation outlined in \S~\ref{subsec:LLS}, the linear regression of the sample $\hat\alpha$ to the fit $(\hat q,\hat m)^T$ corresponds to a linear least--squares problem whose solution is given by the normal equation 
\begin{equation}\label{eq:ordinary_regression}
  \begin{pmatrix}
       \hat q \\
       \hat m
     \end{pmatrix} = \Bigl(B^TB\Bigr)^{-1}B^T\hat\alpha,
\end{equation}
where $B = (1, t_j)_{j=1,...,N}$ is the model matrix and $\hat\alpha_j=\hat\alpha(t_j)$ is the estimated return rate at time $t_j$.
To be an efficient estimator of $(\hat q,\hat m)^T$ the solution \eqref{eq:ordinary_regression} requires that the uncertainty in the input data is homogeneous (i.e. assuming \emph{homoskedasticity}) .
This however is clearly not the case for the return rate sample $\hat\alpha$ as evidenced in Figure~\ref{fig:bias_estimators}(b1--b3); in particular, all three estimators \eqref{eq:OUP}, \eqref{eq:AC1} and \eqref{eq:LLS} appear to have decreasing uncertainty as the system is driven closer to the saddle-node bifurcation.
Therefore the heterogeneity in the uncertainty of the sample $\hat\alpha(t)$ (i.e. \emph{heteroskedasticiy}) makes the ordinary least--square regression \eqref{eq:ordinary_regression} an inefficient point estimate of the linear trend $\hat m$ and intercept $\hat q$ \cite{White1980}.
Assuming that, albeit heterogeneous, the uncertainty $u_{\hat\alpha_j}$ in the estimate $\hat\alpha_j$ is uncorrelated with that of any other estimate $\hat \alpha_k$ in the sequence $\hat\alpha$, then the weighted least--squares formulation
\begin{equation}\label{eq:weighted_regression}
  \begin{pmatrix}
       \hat q \\
       \hat m
     \end{pmatrix} = \Bigl(B^TWB\Bigr)^{-1}B^TW\hat\alpha,\quad W = \mathrm{diag}\Bigl(u^{-2}_{\hat\alpha_j}\Bigr)_{j=1,\dots,N},
\end{equation}
provides the best possible (unbiased) estimator for the linear fit \cite{Aitken1936}.
Using the squared residuals of the ordinary least--squares solution \eqref{eq:ordinary_regression} as heteroskedasticiy--consistent estimators of the entries of the weight matrix $W$ \cite[Thm.~$1$]{White1980}, one can thus implement the full, data-driven (weighted) extrapolation problem \eqref{eq:weighted_regression}, from the sequence $\hat\alpha$, to forecast the breaking time $\hat\tau$.

\section{Practical applications}\label{sec:application}
The tipping time framework and observation protocol, specified in \S~\ref{sec:tipping} and \ref{sec:observation} respectively, provide us with a series of constraints to collect and/or generate timeseries samples in order to apply the proposed, data--driven, forecasting pipeline.
These constraints allow for the optimal design of experiments, data collection and modeling choices for target applications of slow passage through saddle-node B-tipping. 
These specifications are summarized as follows:
\begin{enumerate}
  \item the noisy fluctuations in the residuals $\mathcal{X}$ are additive, i.e. $\sigma(x)=\sigma>0$, and small compared to the timescale separation between the state $x$ and its parameter $\mu$, i.e. $\sigma\ll\sqrt{\varepsilon}$;
  \item the residuals are collected prior to the critical transition and within the admissible time horizon $T_c < \tau - \varepsilon^{-1/3}$;
  \item the residuals are sampled at regular frequency resulting in constant, small enough time intervals $\Delta$ to control bias and uncertainty in the estimation of the return rate $\hat\alpha$;
  \item the residuals are collected over a time window that is long enough so that the correlation in the observation is negligible, i.e. $1/(2\sqrt{\varepsilon(\tau - t)})\ll T_w$, and short enough so that the sample can be considered a stationary approximation of the nonlinear evolution, i.e. $T_w\ll 2(\tau - t)$.
\end{enumerate}
Notice how in the prescribed conditions for the admissible time horizon (ii) and the optimal observation length (iv), prior knowledge of the breaking time $\tau$ is necessary, which is clearly not true in real--world applications.
As a remedy one can deploy an iterative process in which, given a sample, one selects arbitrary values of $T_c$ and $T_w$, estimates the breaking time $\hat\tau$, and then uses such estimate in the next iteration to update the observation protocol accordingly.

In this section we apply the above to design the optimal conditions for which data can be collected to accurately forecast the breaking time of the underlying system.
To do so we consider both low- and high-dimensional cases.
First we apply our pipeline to the canonical dynamic saddle-node normal form \eqref{eq:model} to show how the observation protocol can be used to the design of data collection from a simplified model.
We then subsequentially consider a one--dimensional observable reppresenting the strength of volume transport in the thermohaline ocean circulation generated by a high--dimensional Community Earth System Model (CESM).

\subsection{One--dimensional tipping in the dynamic saddle-node normal form}\label{subsec:dynamic_sn_normal_form}

\begin{figure}[!t]
    \centering 
    \includegraphics[keepaspectratio, width=\textwidth]{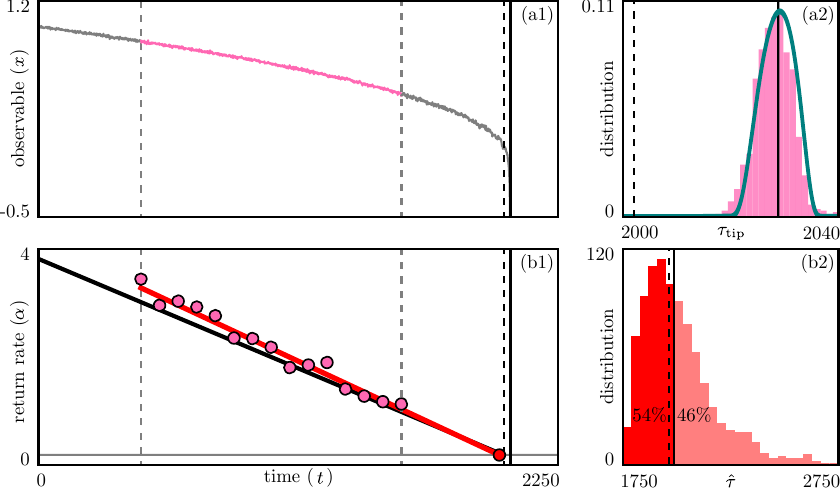}
    \caption{Application of the forecast pipeline to the dynamic saddle--node normal form \eqref{eq:model}. 
      An ensemble $\Omega$ of $N=1000$ sample paths is simulated from $\mu_0=-1$ with additive noise level $\sigma=10^{-2}$.
      Each solution, of which a reppresentative trajectory is show by the grey curve in panel (a1), is slowly driven through the saddle--node bifurcation at time $\tau=2000$, indicated by the black dashed vertical line in panels (a1) and (b1)-(b2), with ramp speed $\varepsilon=5\cdot10^{-4}$.
      The tipping distribution, shown by the pink histogram in panel (a2), follows the small noise approximation $p_1(t)$ given in \eqref{eq:small_noise_approximation}, shown by the teal curve in panel (a2), and its mean approximates the deterministic tipping time $\tau_{\mathrm{tip}}$ of the pullback attractor, depicted by the black solid vertical line in all panels.
      The application of the forecasting pipeline is then applied to the pink subseries in panel (a1), which is extracted from the sample path within the time window $t\in[400,1600]$, depicted by the grey dashed vertical lines in panel (b1), according to the admissible time horizon \eqref{eq:admissible_time_horizon}.
      The estimated return rates $\hat\alpha$, depicted by the pink circles in panel (b1), are computed across strides of a rolling window by solving the ordinary least--squares regression in \eqref{eq:ordinary_least_squares}.
      The forecast of the breaking time $\hat\tau$, reppresented as the red circle in panel (b1), is then obtained by extrapolation, solving the weighted least--square problem \eqref{eq:weighted_regression} as shown by the red linear fit in panel (b1).
      The empirical distribution of forecasted breaking time is reported in panel (b2) with approximately $54\%$ of its mass predicting a breaking time earlier than the deterministic tipping time, highlighted in bright red.
    }
    \label{fig:dynamic_saddle_node}
\end{figure}

We first consider the low--dimensional case of the dynamic saddle--node normal form perturbed by additive noise as in \eqref{eq:model}.
To show how the observation protocol, discussed in \S~\ref{sec:observation}, can be exploited to maximize the forecasting power of the proposed pipeline, we perform numerical simulations that guarantee the slow passage through the bifurcation and delayed tipping in the small noise regime as shown in Figure~\ref{fig:tipping_times_distribution}(a2--b2).

The results of this simulation are depicted in Figure~\ref{fig:dynamic_saddle_node}.
We fix the initial condition at the stable branch $x_+(\mu_0)$ of the equilibria of $f$ for the initial value $\mu_0=-1$ of the bifurcation parameter.
The parameter speed is set at $\varepsilon=5\cdot10^{-4}$ and, from \eqref{eq:breaking_time}, the ground truth of the breaking time is $\tau=2000$ (black dashed vertical line in (a1)).
The tipping time of the pullback attractor (black solid vertical line in (a1)) reads $\tau_{\mathrm{pb}}\approx\tau + 2.338 \varepsilon^{-1/3}\approx 2029.5$, which identifies a short window of opportunity $\Delta\tau\approx29.5$.
The noise intensity is set to $\sigma=10^{-2}$ which satisfies the small noise regime approximation \eqref{eq:small_noise}; this in turn guarantees that the tipping time distribution of an ensemble of simulated sample paths (pink histogram in (a2)) is approximately Gaussian \eqref{eq:small_noise_approximation} (teal curve in (a2)) and the tipping time random variable $\tau_{\mathrm{tip}}$ attains values almost always larger than breaking time $\tau$.

We simulate an ensemble $\Omega = \{\mathcal{X}_k\}_{k=1,\dots,N}$ of $N=1000$ independent sample paths (the grey curve in panel (a1) shows a represenative trajectory); each realization $X_{t_j}^{(k)}\in\mathcal{X}_k$ is sampled at regular frequency $\Delta=10^{-2}$.
According to the observation protocol \S~\ref{subsec:estimators}, from each sample path we extract a subseries (pink curve in panel (a1)) corresponding to the time window $t\in[400,1600]$ (grey dashed vertical lines in panel (a1)), where the lower bound guarantees the extinction of the transient of the initial condition and the upper bound is given by the admissible time horizon \eqref{eq:admissible_time_horizon}.
Each subseries is then analysed in a rolling window to estimate the return rate; from \eqref{eq:window_length} we set the window length to be $T_w = 80$ which correpsonds to $n=8000$ observations for each sample.

For each stride of the rolling window we then use the sample to estimate the return rate \eqref{eq:return_rate} (black curve in panel (b1)) of the underlying, quasi--stationary residuals, which we obtain via linear detrending.
Using the pipeline derived in \S~\ref{subsec:LLS}, we apply the Euler--Maruyama approximation and solve the ordinary least--squares problem \eqref{eq:ordinary_least_squares} obtaining the least biased estimator (see Eq.~\eqref{eq:LLS_bias} and Figure~\ref{fig:bias_estimators}) of the three methods we discussed above.
After a number of strides we obtain a sequence $\hat\alpha$ (pink circles in panel (b1)) of estimated return rates of varying precision (following heteroskedasticity of the estimates).
By \S~\ref{subsec:extrapolators}, we fit a linear model (red curve in panel (b1)) to $\hat\alpha$ by solving the weighted least--squares problem \eqref{eq:weighted_regression} which finally yields the forecast $\hat\tau$ (red circle in panel (b1)).

By repeating the same procedure for any sample path $\mathcal{X}_k$ in the ensemble $\Omega$ we can construct an empirical distribution of the estimated breaking time $\hat\tau$, which is shown by the red histogram in panel (b2).
In this distribution, we highlight those estimates that forecast a breaking time $\hat\tau$ smaller than the tipping time of the pullback attractor $\tau_{\mathrm{pb}}$, which approximates the mean of the ground truth distribution $p_1(t)$ given in \eqref{eq:small_noise_approximation}, indicated by the black solid vertical line.
The total mass of this distribution restricted to $(-\infty,\tau_{\mathrm{pb}})$ is approximately $54\%$ of the total number of estimated breaking times.
Furthermore, we observe that the ensemble mean of the empirical distribution of the estimated breaking time $\E[\hat\tau] = 2057$ is relatively close to the mean of the tipping time distribution $\E[\tau_{\mathrm{tip}}] = 2028$, and $18\%$ of the estimated breaking times fall within $10\%$ of the determinstic breaking time, i.e. within $[\tau_{\mathrm{pb}} - 0.1\Delta T, \tau_{\mathrm{pb}} + 0.1\Delta T]$, where $\Delta T = \tau_{\mathrm{pb}} - T_c$ is the time distance between the admissible time horizon $T_c = 1600$ and the tipping time of the pullback attractor $\tau_{\mathrm{pb}}=2029$.

This entails that, within the simplest model, following an optimal protocol for the generation of data and minimizing the bias in the estimators, the pipeline can still forecast breaking times outside the window of opportunity $\Delta\tau$ with large probability, which might result in delayed or missed warnings.

\subsection{Collapse of the overturning circulation stength in the Atlantic Ocean}\label{subsec:amoc}

The Atlantic Meridional Overturnic Circulation (AMOC) is a volume transport system driving the deep and surface water circulation in the Atlantic Ocean. 
Due to the observed, sustained weakening of its strength, such circulation system has been a major focus of studies in climate science and oceanography communities for the past two decades.
It is hypothesized that the collapse of the AMOC, due to instability-induced transitions, could drastically alter the global climate with particularly catastrophic effects in North--Western Europe \cite{Dijkstra2026}.

The mathematical modeling of the AMOC weakening stength has followed both a low-fidelity, conceptual approach, in which the fundamental dynamics are captured by simplified box--models, and a high-fidelity, computational approach, in which large--scale simulations discretise the North--Atlantic ocean in a relatively fine grid and solve for the governing equations of motions.
Among the former group, we highlight the canonical $2$--box model proposed by Stommel in 1961 \cite{Stommel1961}, the more recent $5$--box extension introduced by Wood et al. in 2019 \cite{Wood2019} and the convective, feedback--delayed multidecadal model of Falkena et al. in 2021 \cite{Falkena2021}; saddle-node B-tipping and R-tipping have been found in all of the above mentioned models \cite{Alkhayuon2019,Lohmann2025a,Mancini2026}.
On the other hand, large--scale, multiphysics simulations, such as those provided by Global Climate Models \cite{vanWesten2023}, are usually too complex to analyze in a systematic and rigorous way; nonetheless saddle-node B-tipping has been hypothesized \cite{vanWesten2025} and early-warning signals have been proposed \cite{vanWesten2024,Smolders2025} providing supporting evidence to the conceptual models that the AMOC is indeed currently on course toward a collapsed state.

\begin{figure}[!t]
    \centering 
    \includegraphics[keepaspectratio, width=\textwidth]{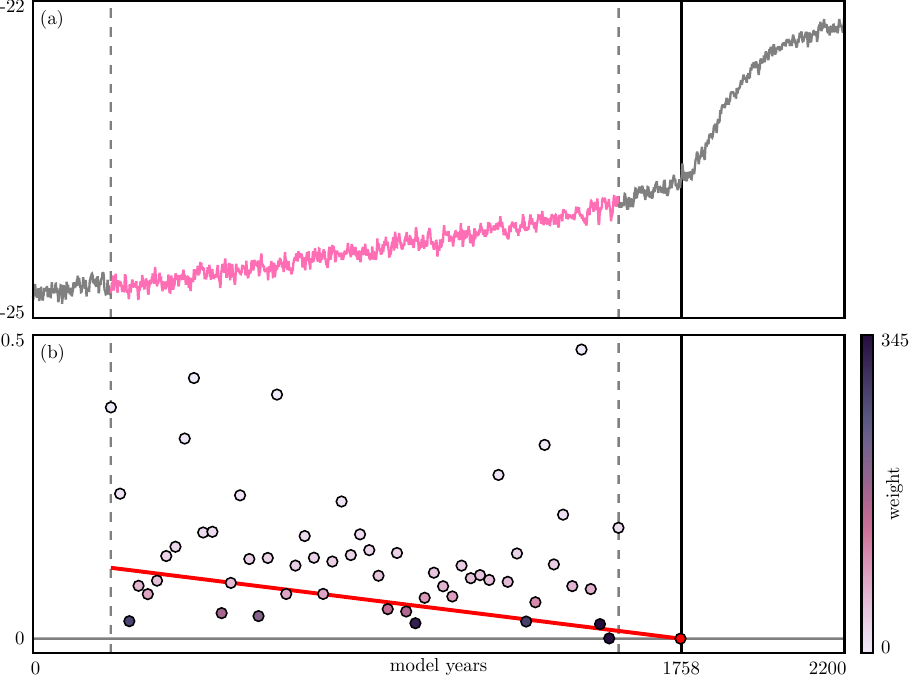}
    \caption{Forecast of the collapse in the overturning circulation strength of the AMOC based on high--dimensional, large--scale simulation in the CESM.
      A univariate observable, depicted by the grey timeseries in panel (a), is computed by subtracting the depth and longitute averages salinity from the depth and longitude averaged temperature at the $34$°S parallel (see \cite{vanWesten2024} for details).
      The observability window, denoted by the vertical dashed grey lines, provides a sample subseries, denoted by the pink timeseries in panel (a), on which the forecast pipeline is applied to.
      The estimated return rates, denoted by the colored circles (colored by their weight in the least--squares problem) in panel (b), are thus obtained by solving the ordinary least--squares problem \eqref{eq:ordinary_least_squares} on rolling windows of the subseries of length $n=25$.
      The solution of the weighted least--squares problem \eqref{eq:weighted_regression}, indicated by the red line in panel (b), is then extrapolated to its zero--crossing yielding a forecasted breaking time at model year $\hat\tau=1755$ (red circle) which anticipates the AMOC collapse by $3$ model years (vertical black line).}
    \label{fig:amoc_collapse}
\end{figure}

In this section we apply the data--driven pipeline introduced in this work to predict the breaking time of the AMOC using large--scale simulations generated by a targeted run of the Community Earth System Model (CESM, version \texttt{1.0.5}) used in phase $5$ of the Coupled Model Intercomparison Project (CMIP5).
Specifically, we use the model's simulation performed in \cite{vanWesten2024} in which a high--dimensional, quasi--steady equilibrium, starting at pre--industrial conditions at model year $0$, is slowly ramped under increasing freshwater forcing in the North--Atlantic. 
The drift of the quasi--steady equilibrium corresponds to the weakining of the circulation strength in the AMOC and it eventually drives the system towards a state of collapse at around model year $1758$.
In \cite{vanWesten2024}, the authours identify a physically--informed early--warning signal of the AMOC collapse based on the freshwater transport measured at the $34$°S.
This has been linked to saddle-node tipping in conceputal climate models \cite{vanWesten2025} and its increase in variability prior to the collapse suggests the possible onset of critical slowing down.
Later studies have confirmed that the $34$°S latitude provides an optimal observation location for the dection of early--warning signals for the AMOC collapse based on the advection of salinity \cite{Smolders2025}.

Mathematically, the freshwater transport is a univariate observable computed from a spatially aggregated, high--dimensional output consisting of the product of the meridional velocity and salinity field averaged across the ocean depth and longitude \cite[Eq.~$(1$b$)$]{vanWesten2024}.
Here we follow a similar approach and compute a univariate observable, namely the circulation strength, from the averaged temperature and salinity timeseries at the proposed latitude.
The circulation strength undergoes a saddle--node bifurcation in Stommel's $2$--box model and interpretable early--warning signals based on such observable have been recently obtained \cite{Papapicco2026}.

\begin{figure}[!t]
    \centering 
    \includegraphics[keepaspectratio, width=\textwidth]{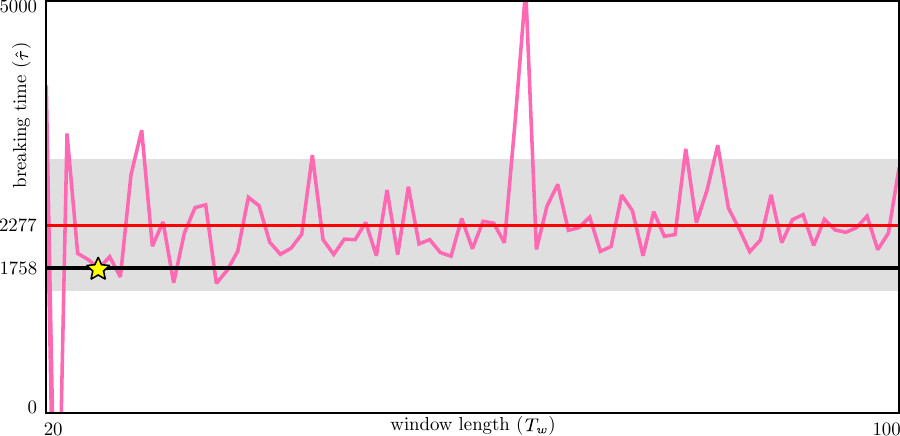}
    \caption{Sensitivity analysis the forecasted breaking time for the AMOC timeseries (see Figure~\ref{fig:amoc_collapse}) upon varying length of the observation window, within the optimal bounds derived in \eqref{eq:window_length} and the admissible time horizon \eqref{eq:admissible_time_horizon}.
      The ground truth of the model $\tau=1758$ is reported by the black line while the forecast $\hat\tau=1755$ for the optimal window length $T_w = 25$, used in Figure~\ref{fig:amoc_collapse} is depicted by the yellow star.
      The mean over the varying window lengths is depicted by the red line while the realizations bounded within one standard deviation away are seen within the grey band.}
    \label{fig:sensitivity_analysis}
\end{figure}

The result of our experiment are depicted in Figure~\ref{fig:amoc_collapse}.
The circulation strength solution, indicated by the grey timeseries in panel (a) resembles the drift of a quasi--steady equilibrium driven through the saddle--node bifurcation which, following \cite{vanWesten2024}, is placed at model year $\tau=1758$, indicated by the vertical black lines.
In order to extract a reasonable signal for our forecast problem, we follow the observation protocol discussed in \S~\ref{sec:observation}. 
In the model run, the freshwater forcing is increased linearly at a rate $\varepsilon=3\cdot10^{-4}$ and solutions are recorded at a frequency of $\Delta=1$ year.
This entails that the admissible time horizon can be safely placed at $T_c=1600$, depicted by the dashed vertical grey lines, while the rolling window's length is set to be $T_w=25$ model years (corresponding to $n=25$ observations).
As for the previous example, at each stride of the window we solve the ordinary least--squares problem \eqref{eq:ordinary_least_squares} to estimate the return rate, resulting in a sequence $\hat\alpha$ depicted by the pink circles in panel (b).
Solving the weighted least--squares regression \eqref{eq:weighted_regression} we thus fit a linear model to the estimates $\hat\alpha$, reppresented by the red curve, whose slope and intercept yield a forecasted breaking time, depicted by the red circle.
Our forecast, which predicts the AMOC collapse to be at model year $\hat\tau=1755$, succesfully anticipates the tipping time of the simulated large--scale model, providing a window of opportunity of $\Delta\tau=3$ model years.

While this result seemingly strengthens the computational evidence supported in the literature for the same model run \cite{vanWesten2024,Smolders2025}, care shall be taken in its interpretation.
Given the current values of $\varepsilon$ and $\tau$ in fact, we have access to a large interval of optimal window lengths to choose for the forecasting problem.
In Figure~\ref{fig:amoc_collapse} we choose the fixed width $T_w=25$ which is relatively closer to the lower bound derived in \eqref{eq:lower_bound}.
By varying $T_w\in[20,100]$ and applying our extrapolation pipeline repeatedly, we can analyse the sensitivity of the forecast to the free parameter that consitutes the window length.
The results of this analysis are depicted in Figure~\ref{fig:sensitivity_analysis}.
The mean forecasted breaking time (red line) is significantly larger (approximately $2277$ model years) than the optimally forecasted time (yellow star) for the window length $T_w=25$ reported in Figure~\ref{fig:amoc_collapse}.
Nevertheless approximately $88\%$ of those windows length predict a forecasting time that is within one standard deviation away from such mean (grey band).

\section{Discussion and conclusions}\label{sec:conclusion}
In this manuscript we have derived a fully data--driven pipeline for the forecast of the time of collapse in complex systems, commonly known as critical transitions or tipping events.
Our pipeline estimates the return rate of the quasi--stationary residuals of a windowed timeseries and then extrapolates the estimates sequence to its zero--crossing by fitting a linear model.

We have shown that a quasi--maximum likelihood estimation of the Euler--Maruyama approximation of the recorded residuals, resulting in the solution of an ordinary least--squares problem, results in a more accurate and precise estimation of the return rate compared to the more established sample variance \cite{Kuehn2011} and autocorrelation at lag--$1$ \cite{Boers2021}.
Furthermore, we derive bounds on the optimal length for the sample to yield meaningful estimations under the stationary Ornstein--Uhlenbeck assumption.

Monte Carlo simulations show that the proposed technique predicts the tipping time of a sample, by extrapolating the linear decay of the estimated return rates, before the actual collapse of the system in more than half of the cases.
Furthermore, we were able to estimate the tipping time with great accuracy in approximately one-fifth of the trials.
This is further evidenced by the application of the forecast pipeline to a large--scale simulation of the AMOC collapse, which we predict with $3$ model years in advance before the ground truth.
Sensitivity analysis however showed that results are widely affected by the changes in the window length despite it being bounded by the optimal protocols we derive in the paper.

In the construction of the methodology we have leveraged results in nonautonomous dynamical systems theory and well-established techniques in the statistical inference of stochastic processes.
As a result, the succesful application of this pipeline to timeseries samples relies on a set of assumptions on the quality of the data and the underlying nature of the systems being observed.
First, we expect the sample to be an observable of a high--dimensional complex system, that is slowly driven through a saddle--node bifurcation and that is perturbed by additive noise, to be optimal in the sense prescribed by \cite{Lohmann2025b,Ditlevsen2026}.
Second, we expect that the level of noise in the observations has relatively small magnitude compared to the speed at which the system is ramped toward the catastrophe, so to minimize the likelihood of early, noise-induced escapes from the basin of attraction.
Third, the sample is recorded at regular frequency and it is sufficiently distant from the saddle--node bifurcation so that the statistical estimation of the return rate is feasible.

So far the theory of early-warning signals, including the proposed pipeline in this manuscript, has largely been developed within the same setup outlined above, which leaves a large set of models for which theoretical develpment is lacking.
The relaxation of the above assumptions thus requires extensive mathematical theory; recent early development being made in this direction include works on systems driven by multiplicative noise \cite{Morr2024b} and jump processes \cite{Layritz2025}, nonlinear \cite{Ritchie2016} and periodically forced \cite{Kuehn2018} parameter shifts with little or no timescale separation, and spatially extended models \cite{Gowda2015,Bernuzzi2024,Clarke2026,Miyauchi2026}.

\printbibliography

\end{document}